\documentclass{article}
\usepackage{amsmath}
\usepackage{amssymb}
\usepackage{amsthm}
\usepackage{fullpage}

\newtheorem{Theorem}{Theorem}
\newtheorem{Proposition}{Proposition}
\newtheorem{Lemma}{Lemma}
\newtheorem{Corollary}{Corollary}
\newcommand{\convd}{\mbox{$ \ \stackrel{\!d}{\rightarrow}$ }}
\newcommand{\eqd}{\mbox{$ \ \stackrel{\!d}{=}$ }}
\newcommand{\be}{\begin{equation}}
\newcommand{\ee}{\end{equation}}
\newcommand{\ba}{\begin{eqnarray}}
\newcommand{\ea}{\end{eqnarray}}

\def\P{\mathbb{P}}
\def\z{\zeta}

\def\s{\sigma}
\def\L{\mathcal{L}}
\def\lf{\lfloor}
\def\rf{\rfloor}

\begin{document}

\title{Colored loop-erased random walk on the complete graph
\thanks{%
Research supported in part N.S.F.\
grant DMS-0405779}}
\author{Jomy Alappattu
\thanks {Dept. Mathematics, U. C. Berkeley; email jomy@math.berkeley.edu} \and %
Jim Pitman
\thanks{%
Dept. Statistics, U. C. Berkeley; email pitman@stat.berkeley.edu}
}
\maketitle

\begin{abstract}
Starting from a sequence regarded as a walk through some set of
values, we consider the associated loop-erased walk as a sequence of
directed edges, with an edge from $i$ to $j$ if the loop erased walk
makes a step from $i$ to $j$.  We introduce a coloring of these
edges by painting edges with a fixed color as long as the walk does
not loop back on itself, then switching to a new color whenever a
loop is erased, with each new color distinct from all previous
colors. The pattern of colors along the edges of the loop-erased
walk then displays stretches of consecutive steps of the walk left
untouched by the loop-erasure process. Assuming that the underlying
sequence generating the loop-erased walk is a sequence of
independent random variables, each uniform on $[N]:=\{1, 2, \ldots,
N\}$, we condition the walk to start at $N$ and stop the walk when
it first reaches the subset $[k]$, for some $1 \leq k \leq N-1$.  We
relate the distribution of the random length of this loop-erased
walk to the distribution of the length of the first loop of the
walk, via Cayley's enumerations of trees, and via Wilson's
algorithm. For fixed $N$ and $k$, and $i = 1,2, \ldots$, let $B_i$
denote the event that the loop-erased walk from $N$ to $[k]$ has $i
+1$ or more edges, and the $i^{th}$ and $(i+1)^{th}$ of these edges
are colored differently. We show that given that the loop-erased
random walk has $j$ edges for some $1\leq j \leq N-k$, the events
$B_i$ for $1 \leq i \leq j-1$ are independent, with the probability
of $B_i$ equal to $1/(k+i+1)$. This determines the distribution of
the sequence of random lengths of differently colored segments of
the loop-erased walk, and yields asymptotic descriptions of these
random lengths as $N \to \infty$.
\end{abstract}

\section{Introduction}
The {\it loop-erased walk} derived from a sequence $(X_n, n =
0,1,\ldots)$ is a sequence $(\L_n, n = 0,1,\ldots)$ of finite
subsequences of $(X_n, n = 0,1,\ldots)$ defined as follows: let
$\L_0 = (Y_{0,0}) = (X_0)$, and inductively, if $X_n$ is not in
$\L_{n-1} = (Y_{n-1,0}, \ldots, Y_{n-1,L_{n-1}})$, then form $\L_n$
by appending $X_n$ to the end of $\L_{n-1}$, i.e., $\L_n = (Y_{n,0},
\ldots, Y_{n,L_n})$ with $L_n = L_{n-1} + 1$, $Y_{n,i} = Y_{n-1,i}$
for $0 \leq i \leq L_{n-1}$, and $Y_{n,L_n} = X_n$. On the other
hand, if $X_n = Y_{n-1,j}$ for some $1 \leq j \leq L_{n-1}$, then
construct $\L_n$ by truncating the part of $\L_{n-1}$ beyond
$Y_{n-1,j}$, i.e., let $L_n = j$, and define $\L_n = (Y_{n,0},
\ldots, Y_{n,L_n})$ by $Y_{n,i} = Y_{n-1,i}$ for $0 \leq i \leq j$;
in this case, we say that a {\it loop} has occurred. For each $n =
0, 1, \ldots$, we call $L_n$ the {\it length of the loop-erased walk
at time $n$}, with the understanding that if $\L_n$ is a single
point $(X_i)$ for some $i$, the length is zero. So the length
measures the number of steps of the path, or the number of {\em
edges} of the path, with each edge representing some one-step
transition $(X_n, X_{n+1})$ of the original sequence
$(X_n,n=0,1,\ldots)$. See \cite{Law} and \cite{LP} for equivalent
alternative definitions of loop-erased walk, and discussions of some
basic results on the loop-erasure of random sequences.

There is a natural way to ``color'' the loop-erasure of the walk as follows.
Assume that we have some
infinite palette $\{C_1, C_2, \ldots \}$ of colors. Run the walk,
and until the first loop occurs, color the edges of the walk with
the color $C_1$. When the first loop occurs, erase the colored edges
as the definition of loop-erasure requires, and continue the walk,
now coloring the subsequent edges with the color $C_2$ until the
next loop occurs. More generally, keep coloring the edges of the
walk with a fixed color $C_i$ until a loop occurs, at which point we
change to a new fixed color $C_{i+1}$.

Fix a positive integer $N\geq 2$, and let $(X_n, n =0,1,\ldots)$ be
an sequence of independent random variables, with $X_0 = N$ for
convenience, and $X_1, X_2, \ldots$ independent and uniformly
distributed on the set $[N]:= \{1,2,\ldots,N\}$. We use this random
sequence to construct a {\it loop-erased random walk} on $[N]$,
following the definition for loop-erased walk above. Apart from some
delay due to self-loops when $X_{n+1} = X_n$ for some $n$, the
sequence of steps of the loop-erased walk is the same as if it were
derived from a random walk on the complete graph. In particular, we
are interested in the random coloring of the loop-erased walk when
it first reaches the subset $[k]$ for some $1 \leq k \leq N-1$, and
this random coloring has the same distribution for a sequence of
independent random variables $(X_n, n=0,1,\ldots)$ as for a random
walk on the complete graph.
Let $R_N$ denote the {\it first repeat time} for $X_0, X_1, \ldots$,
i.e., the first index $i$ such that $X_i \in \{X_0, \ldots,
X_{i-1}\}$, which is the length of the first loop that is erased in
the process of loop-erasure. The distribution of $R_N$ is determined
by the well-known solution of the classical birthday problem, that
is,
\ba \label{rng} \P(R_N > j) &=& \prod_{i=1}^j \frac{N-i}{N}
\\ \label{rne} \P(R_N = j) &=& \frac{j}{N} \prod_{i=1}^{j-1}
\frac{N-i}{N}.\ea
Our main result relates the distribution of $R_N$
to the distribution of the length of the loop-erased walk stopped
when it hits the target set $[k]$.

Let $\zeta_{N,k}$ be the first time $t$ such that $X_t \in [k]$, and
note that $\zeta_{N,k}$ is a geometric random variable with
parameter $k/N$. We use the notation $X \eqd Y$ to mean that random
variable $X$ and $Y$ have the same distribution, and $X \eqd (Y|A)$
to indicate that the distribution of $X$ is the same as the
conditional distribution of $Y$ given the event $A$. Also let
$$(a)_b := a(a-1)(a-2) \cdots (a-b+1)$$
be the usual falling factorial for $b = 1, 2,\ldots$ with $(a)_0 =
1$.

\begin{Theorem} \label{seg}
Let $\lambda_{N,k}$ be the length of the
loop-erased random walk derived from $X_0 = N$ and a sequence of independent variables $X_1, X_2, \ldots$
with uniform distribution on $[N]$, stopped at time $\zeta_{N,k}$ when the sequence first hits $[k]$,
and let $R_N$ be the first repeat time derived from the same sequence.
Then
$$\lambda_{N,k} \eqd (R_N - k |R_N > k);$$
 that is, for every $1 \leq j \leq N-k$,
\begin{equation} \label{lrnk}
\P(\lambda_{N,k} = j) = \P(R_N - k = j| R_N > k) = \frac{(k+j)(N-k-1)_{j-1}}{N^j}
\end{equation}
Moreover, in the colored loop erased walk of length $\lambda_{N,k}$
obtained by stopping at time $\zeta_{N,k}$, let $B_i$ denote the
event that the $i^{th}$ and $(i+1)^{th}$ edges of the loop-erased
walk are colored differently. Then given $\lambda_{N,k} = j$, the
events $B_i$ for $1 \leq i \leq j-1$ are independent with
$$
\P(B_i) = \frac{1}{k+i+1}.
$$
\end{Theorem}

We prove this result in Section \ref{prf}, then show in Section
\ref{ust} how the simple formula (\ref{lrnk}) is closely related
both to Wilson's loop-erased random walk algorithm to generate
spanning trees of a graph, and to Cayley's formula for the number of
forests with a fixed number of vertices and a fixed set of roots.
In Section \ref{scalims}, we relate Theorem \ref{seg} to other basic
results about compositions which are closely connected to Aldous's
Brownian continuum random tree \cite{Al1} and to stick-breaking
schemes. In Section \ref{open}, we discuss some open questions which
arise naturally from our analysis.

\section{Proof of Theorem \ref{seg}} \label{prf}

Most of the work for this proof is done by the following lemma,
where we use the notation of Theorem \ref{seg}:

\begin{Lemma}\label{switches}

Let $\nu_{N,k}$ denote the number of different colors of segments of
the colored loop-erased walk started at $X_0 = N$ and stopped on
first reaching $[k]$, let $\lambda_{N,k}$ be the total length of the
loop-erased walk at this time, and given $\lambda_{N,k} = j$, let
$Y_i$ denote the indicator of the event $B_i$ that the $i^{th}$ and
$(i+1)^{th}$ edges in the loop-erased random walk from $N$ to $[k]$
are colored differently. Then for each choice of positive integers
$s$ and $j$ with $s \le j$ and  each choice of positive integers
$i_1, \ldots, i_{s-1}$ with
$$
0 < i_1 < \cdots < i_{s-1} < j
$$
\ba \nonumber &\P(\lambda_{N,k} = j, \nu_{N,k} = s, Y_{i_1} = 1, Y_{i_2}
= 1, \ldots, Y_{i_{s-1}} = 1, Y_i = 0 \text{ for } i \in [j-1] - \{i_1, \ldots, i_{s-1}\})\\
\label{twolem}=& \frac{(k+1)(N-k-1)_{j-1}}{N^j} \prod_{r=1}^{s-1}
\frac{1}{k+i_r} \ea where if $s = 1$, the product is understood to equal 1.
\end{Lemma}

\begin{proof}
Consider first the case when all edges are painted the same color,
i.e., $\nu_{N,k} = s = 1$, and $Y_i = 0$ for all $i$. Then
loop-erasure can only occur at vertex $N$, and thus, the random walk
takes some number $n$ of steps to vertices in $[N]-[k]$, then it
hits $N$ one final time, and then it hits $(j-1)$ distinct vertices
excluding $N$ and $[k]$ before finally hitting some element of
$[k]$. It thus follows that that \begin{eqnarray*} \P(\lambda_{N,k}
= j, \nu_{N,k} = 1) &=& \left(1 + \sum_{n=0}^\infty \left(1 -
\frac{k}{N}
\right)^n\frac{1}{N}\right) \frac{(N-k-1)_{j-1}}{N^{j-1}} \frac{k}{N} \\
&=& \frac{(k+1)(N-k-1)_{j-1}}{N^{j}}\end{eqnarray*} as desired.

Next, consider the case of two colors, with a color change at the
$h^{th}$ vertex in the path from $N$ to $[k]$ for some $1 \leq h
\leq j-1$. Again, the random walk takes $n$ steps outside of $[k]$
before hitting vertex $N$ for the last time, and then the random
walk takes $h$ steps without looping, to distinct vertices. At the
$h^{th}$ vertex, the random walk takes $h$ steps to vertices outside
of $[k]$ and the $0^{th}, 1^{st}, \ldots, (h-1)^{th}$ vertices
before hitting the $h^{th}$ vertex one final time, and then it takes
$j-h$ steps without looping until it hits something in $[k]$. So
again appealing to independence and the uniform distribution of the
variables $X_i$ for $i \ge 1$, we see that \begin{eqnarray*}
&\P(\lambda_{N,k} = j, \nu_{N,k} = 2, Y_h = 1, Y_i = 0 \text{ for }
i \in [j-1] - \{h\})\\=&  \left(1 + \sum_{n=0}^\infty \left(1 -
\frac{k}{N} \right)^n \frac{1}{N} \right) \left(\sum_{m=0}^\infty
\left(1 - \frac{k+h}{N}\right)^m \frac{1}{N} \right)
\frac{(N-k-1)_{j-1}}{N^{j-1}} \frac{k}{N} \\=&
\frac{(k+1)(N-k-1)_{j-1}}{N^j} \cdot \frac{1}{k+h} \end{eqnarray*}
as desired.

Extending this argument to three or more colors is straightforward.
\end{proof}

\begin{proof}[Proof of formula (\ref{lrnk})]
The equality of the second and third
expressions is evident from (\ref{rng}) and (\ref{rne}).
To prove
the equality of the first and third expressions in this equation,
we sum up equation (\ref{twolem}) over all $s$ between
1 and $N-k$ and all possible sequences $(i_1, \ldots, i_{s-1})$,
corresponding to all possible subsets of $[j-1]$, including the empty subset, and to
all possible sequences of values of the color change indicators ($Y_i$).
Then (\ref{lrnk}) is seen to amount to
\begin{equation}\label{fj}
f(j) := \sum_{I \subset [j-1]} \prod_{i \in I} \frac{1}{k+i} = \frac{k+j}{k+1}.
\end{equation} where it is understood that if $I$ is the empty set, then the product
equals 1. It is obvious that $f(1) = 1$; suppose inductively that $f(j-1) = \frac{k+j-1}{k+1}$.
Then using the fact that $f(j) = f(j-1) + \frac{1}{k+j-1}f(j-1)$ yields equation (\ref{fj})
for any $j$, and finishes the proof of (\ref{lrnk}).
\end{proof}

For the second part of the Theorem, we need need a lemma about
Bernoulli trials:

\begin{Lemma}\label{bern}
A sequence $Y_1, \ldots, Y_{j-1}$ is an independent sequence of Bernoulli
random variables, such that $Y_i$ has parameter $\frac{1}{k+i+1}$, if and only if
for each choice of positive integers $i_1, \ldots, i_{s-1}$ with
$$
0 < i_1 < \cdots < i_{s-1} < j,
$$
$$P(Y_{i_1}
=\cdots = Y_{i_{s-1}} = 1, Y_i = 0 \text{ for } i \in [j-1] - \{i_1, \ldots, i_{s-1}\})
= \frac{k+1}{k+j}\prod_{r=1}^{s-1} \frac{1}{k+i_r}.$$
\end{Lemma}
\begin{proof}
Suppose that the sequence $Y_1, \ldots, Y_{j-1}$ is independent Bernoulli, such
that $Y_i$ has parameter $p_i:=\frac{1}{k+i+1}$. Then $$P(Y_{i_1}
=\cdots = Y_{i_{s-1}} = 1, Y_i = 0 \text{ for } i \in [j-1] - \{i_1, \ldots, i_{s-1}\})
= \prod_{r=1}^{s-1} p_{i_r} \prod_{i \in [j-1] - \{i_1, \ldots, i_{s-1}\}}
\left(1 - p_i\right)$$ Using the fact that $p_i = p_{i-1}(1-p_i)$ for $i = 1,
\ldots, j-1$ (where we let $p_0 = \frac{1}{k+1}$) and the fact that $(1-p_1) \cdots (1-p_{j-1})
= \frac{k+1}{k+j}$, we see that this last product becomes $$\prod_{r=1}^{s-1} p_{i_r - 1}
\prod_{i=1}^{j-1} (1-p_i) = \frac{k+1}{k+j} \prod_{r=1}^{s-1} p_{i_r-1} = \frac{k+1}{k+j}
\prod_{r=1}^{s-1} \frac{1}{k+i_r},$$ as desired. The converse is obvious by just
reversing the sequence of equalities.
\end{proof}

\begin{proof}[Proof of Theorem \ref{seg}]
Formula (\ref{lrnk}) has already been established.
From Lemma \ref{switches} and (\ref{lrnk}) we see that \ba \nonumber &\P(\nu_{N,k} = s,
Y_{i_1} = 1, Y_{i_2} = 1, \ldots, Y_{i_{s-1}} = 1, Y_i = 0
\text{ for } i \in [j-1] - \{i_1, \ldots, i_{s-1}\}|\lambda_N = j)\\&=\frac{k+1}
{k+j} \prod_{r=1}^{s-1} \frac{1}{k+i_r} \ea and now the
conclusion follows from Lemma \ref{bern}.
\end{proof}

\section{The length of the loop-erased random walk} \label{ust}
\subsection{The Markov property of the length of the loop-erasure}
An alternative method of proving formula (\ref{lrnk})
begins with
the following observation: if $L_n$ denotes the length of the loop-erasure of the i.i.d. sequence $(X_0, \ldots, X_n)$
of random variables uniform on $[N]$,
then $(L_n, n = 0,1,\ldots)$ has the same dynamics as a Markov chain with
the following transition probabilities, started at $L_0 = 0$: \be \label{qN} Q_N(i,j) =
\left\{\begin{array}{ll} 1/N & 0 \leq j \leq i
\\ (N-i-1)/N & j = i + 1 \\ 0 & \text{otherwise} \end{array}\right.\ee
In fact, for $n = 1,2,\ldots$, $L_n \eqd \min(L_{n-1}+1, X_n-1)$.
Although by definition $L_0 = 0$ throughout this article, we could just as well
start with some arbitrary loop-less path (and a vertex at the end of the path from which to step)
whose length $L_0$ is a random variable taking values in $\{0, 1, \ldots, N-1\}$
and then run loop-erased random walk; if we
then choose $X_0$ independent of $X_1, X_2, \ldots$ so that
$X_0 -1\eqd L_0$, it follows by induction that $$L_n \eqd \min_{0 \leq j \leq n} (X_{n-j} + j-1).$$
Using the independence of the $(X_i, i = 0,1,\ldots)$, it follows that
\ba \nonumber Q^n_N(i,[m,\infty)) &=& \P(L_n \geq m|L_0 = i) \\ \nonumber
&=& \P(\min_{0 \leq j \leq n}(X_{n-j} + j-1) \geq m|X_0 = i) \\ &=& \label{qnn}
1(i+n -1\geq m) f_{N,m-n+1} f_{N,m-n+2} \cdots f_{N,m} \ea where for $X$ a
random variable uniformly distributed on $[N]$, $f_{N,m} := \P(X > m)$;
note that $f_{N,m} = 1$ if $m \leq 0$ and $f_{N,m} = 0$ if $m \geq N$.
From this, we obtain arbitrary entries of powers of the transition matrix:
\begin{eqnarray}
\nonumber Q^n_N(i,m)&=& Q^n_N(i,[m,\infty)) - Q^n_N(i,[m+1,\infty))\\&=& \label{qnim}
f_{N,m-n+2} \cdots f_{N,m} (1(i + n-1 \geq m) f_{N,m-n+1} - 1(i+n-1 \geq m+1) f_{N,m+1})
\end{eqnarray}
Moreover, letting $n \to \infty$ in (\ref{qnn}) and using the fact that $X_1,X_2, \ldots$
are nonnegative random variables, it follows that
\begin{eqnarray*}\lim_{n \to \infty}Q^n(i,[m,\infty)) &=& f_{N,1} \cdots f_{N,m} \\
\lim_{n \to \infty}Q^n(i,m) &=& f_{N,1} \cdots f_{N,m} (1 - f_{N,m+1})\end{eqnarray*}
where if $m \leq 0$, the product $f_{N,1} \cdots f_{N,m}$ is understood to be 1.
Note that the first limit is the probability
$\P(R_N > m)$, and the second limit is the probability
$\P(R_N = m+1)$. Thus, applying the convergence theorem for irreducible aperiodic Markov chains \cite[page 314]{Dur},
we obtain the following result: \begin{Proposition}
The stationary distribution of the Markov chain $(L_n, n = 0,1,\ldots)$ is the distribution
of the random variable $R_N-1$, where $R_N$ is the index of the first repeat in an i.i.d. sequence
of random variables uniform on $[N]$. \end{Proposition}

As an aside, the exact same reasoning can be applied to a non-uniform
random variable $X$ on the positive integers, to obtain the following:
\begin{Corollary}
Let $X$ be a positive-integer-valued random variable, and define
an independent sequence $X_0, X_1, \ldots,$, where $X_0$ has some distribution
on the positive integers, and $X_1, X_2, \ldots$ is an i.i.d. sequence of variables with
the same distribution as $X$. Define a transition matrix $Q$ on the
nonnegative integers as follows:
$$Q(i,j) = \left\{\begin{array}{ll} \P(X = j+1) & 0 \leq j \leq i \\
\P(X > i+1) & j = i+1\end{array}\right.$$ Then if $L_n$
is the loop-erasure of the path $(X_0, \ldots, X_n)$, $(L_n, n = 0,1,\ldots)$
has the same transition dynamics as a Markov chain with
transition matrix $Q$, and if $g_m := \P(X > m)$, then powers of the transition matrix
are given by $$(Q)^n(i,m) = g_{m-n+2} \cdots g_m(1(i+n-1 \geq m)g_{m-n+1} - 1(i+n-1\geq m+1)g_{m+1}).$$
Moreover, if $\P(X=1) > 0$, then the Markov
chain is irreducible and positive recurrent on the nonnegative integers, with a stationary
distribution determined by either of the following formulas for
$m =1,2,\ldots$:
\begin{eqnarray*}\pi([m,\infty))&=& \prod_{i=1}^m \P(X > i)\\
\pi(m) &=& \P(X \leq m+1) \prod_{i=1}^m \P(X > i)
\end{eqnarray*}
\end{Corollary}

We now use equation (\ref{qnim}) to provide a second proof of formula (\ref{lrnk}). Conditioning
on the value of $\zeta_{N,k}$, which is a geometric random variable with
parameter $k/N$, we see that for $1 \leq j \leq N-k$,
\begin{eqnarray}
\nonumber \P(\lambda_{N,k} = j) = \P(L_{\zeta_{N,k}} = j) &=& \sum_{n=1}^\infty \P(L_{n-1} = j-1|\zeta_{N,k} = n)
\P(\zeta_{N,k} = n) \\ \label{pnkj}&=& \sum_{n=1}^\infty Q^{n-1}_{N-k}(0,j-1) \left(1 - \frac{k}{N}\right)^{n-1}\frac{k}{N}
\end{eqnarray}
To calculate $Q^{n-1}_{N-k}(0,j-1)$, we use equation (\ref{qnim}), to see that
\begin{equation}
Q^{n-1}_{N-k}(0,j-1) = \left\{\begin{array}{ll} f_{N-k,1} \cdots f_{N-k,j-1}(1 - f_{N-k,j})
=\P(R_{N-k} = j) &\text{ if }
0 \leq j \leq n-2 \\ f_{N-k,1} \cdots f_{N-k,j-1} = \P(R_{N-k} > j-1) &\text{ if }j = n-1 \\ 0&\text{ if } j \geq n \end{array}\right.
\end{equation}
Now returning to equation (\ref{pnkj}), we see that
\begin{eqnarray}
\P(\lambda_{N,k} = j) &=& \P(R_{N-k} > j-1)\left(1 - \frac{k}{N}\right)^{j-1}\frac{k}{N}
+\P(R_{N-k} = j)\sum_{n=j+2}^\infty \left(1 - \frac{k}{N}\right)^{n-1} \frac{k}{N}\\
&=&\frac{(N-k-1)_{j-1}}{(N-k)^{j-1}}\cdot \frac{(N-k)^{j-1}}{N^{j-1}}\cdot \frac{k}{N} + \frac{j}{N-k}
\cdot \frac{(N-k-1)_{j-1}}{(N-k)^{j-1}} \cdot \frac{(N-k)^j}{N^j} \\&=&
\frac{(j+k)(N-k-1)_{j-1}}{N^j},
\end{eqnarray}
as desired.

\subsection{Relation to using Wilson's algorithm} \label{wilson}
Yet another derivation of formula (\ref{lrnk}) is provided by Wilson's
algorithm \cite{Wil}; to explore this connection, we will need to introduce
some preliminaries on trees. A {\it rooted tree} $T = (V,E,r)$ consists of a vertex
set $V$, an edge set $E \subset V \times V$, and a distinguished vertex $r \in V$
called the {\it root}, such that for any non-root vertex $v \in V$, there
is a unique directed sequence of edges that leads from $v$ to $r$, and such that
there are no undirected loops, i.e., from any vertex there does not exist a sequence
of distinct undirected edges which leads back to that vertex.

Wilson's algorithm can make use of the random walk on the complete graph
to generate a random tree with $N$ vertices labeled by $[N]$ as follows:
Let $T_0$ be the one-point tree with root and vertex labeled $1$. Suppose that $T_0, \ldots, T_{n-1}$
have been defined, with respective vertex sets $V_0 = \{1\}, V_1, \ldots, V_{n-1}$.
If $T_{n-1}$ has $N$ vertices, stop the algorithm
and output $T_{n-1}$. Otherwise, pick some vertex $v \in [N] - V_{n-1}$ (it does not
matter how one chooses $v$), and run a random walk on the complete graph from $v$ until it hits
some vertex in $V_{n-1}$. Loop-erase this random walk, and add the loop-erased path
to $T_{n-1}$ to form the tree $T_n$, still with root labeled $1$, and vertex set $V_n$.
Also, call this loop-erased path from $w$ to $V_{n-1}$ a {\it macrostep}
of the algorithm.

According to \cite{Wil}, the random tree generated by Wilson's
algorithm applied to the complete graph with $N$ vertices is
uniformly distributed among all rooted trees labeled by $[N]$ with
root $1$; call this the {\it uniform spanning tree with root} 1
(where the word {\it spanning} is used to signifiy that the tree has
the full set of $N$ vertices). Suppose that we start Wilson's
algorithm at a vertex $N$. The first macrostep is just the
loop-erased walk from $N$ to $1$, and thus the length of this
macrostep has the same distribution as $\lambda_{N,1}$. On the other
hand, Wilson's algorithm implies that this macrostep is also the
path from $N$ to $1$ in the uniform spanning tree with root 1. If
$H_{N,1}$ is the length of this path, i.e., the number of edges in
the path, then Wilson's algorithm clearly implies that $H_{N,1} \eqd
\lambda_{N,1}$. Moreover, Meir and Moon \cite{MM} proved that
$H_{N,1} \eqd (R_N-1|R_N > 1)$, and thus Wilson's algorithm coupled
with this result in \cite{MM} yields an alternative proof of formula
(\ref{lrnk}) in the case $k = 1$. In fact, the methods of Wilson's
algorithm and Meir and Moon can be applied to a random {\it rooted
forest labeled by $[N]$} (i.e., a collection of trees with $N$ total
vertices labeled by $[N]$) with a fixed set of roots labeled by
$[k]$, which is uniform among all such rooted forests labeled by
$[N]$ with the same root set $[k]$, to prove formula (\ref{lrnk})
for arbitrary $k$.

The result of Meir and Moon made use of Cayley's formula for the enumeration of
forests. Thus our derivation of formula (\ref{lrnk}) yields
a novel proof of Cayley's formula:

\begin{Corollary}[Cayley's formula] The number $t_{N,k}$ of forests with $N$ vertices
labeled by $[N]$ and $k$ rooted trees with root set $[k]$ is given by $t_{N,k} = kN^{N-k-1}$.
\end{Corollary}

\begin{proof}
As mentioned above, applying Wilson's algorithm with root $1$, started at some other vertex
$v$, proves that $\lambda_{N,1}$---the length of the loop-erased random walk
from $v$ to $1$---and $H_{N,1}$---the edge-distance between $v$ and $1$ in a spanning tree
labeled by $[N]$ which is uniform among all spanning trees with root $r$---have the same
distribution. For any fixed $1 \leq j \leq N-1$, we want to count the number
of spanning trees with root $r$ such that the edge-distance from
$v$ to $r$ equals $j$. We have to choose the $j-1$ vertices in the path
from $v$ to $r$, and then each vertex on the path (including $v$ and $r$) may be considered
the root of a tree; thus, there are $(N-2)_{j-1} t_{N,j+1}$ such spanning trees.
Thus, it follows that \be \label{usteq} \frac{(j+1)(N-2)_{j-1}}{N^j} =
\P(\lambda_{N,1} = j) = \P(H_{N,1} = j) = \frac{(N-2)_{j-1} t_{N,j+1}}{t_{N,1}}\ee
In the particular case $j = N-1$, it is obvious that $t_{N,j+1} = 1$, so this equality
implies that $t_{N,1} = N^{N-2}$. But now substituting this into equation (\ref{usteq})
yields Cayley's formula for arbitrary $j$.
\end{proof}

Lyons and Peres \cite{LP} use very similar reasoning in applying Wilson's
algorithm to calculate $t_{N,1}$, by computing the probability of a
particular tree, namely a path of length $N-1$.

\section{Some scaling limit results} \label{scalims}
\subsection{Connections to Poisson and Rayleigh processes}
See Pittel's paper \cite{Pit} for a discussion of related results
where a similar distribution involving the lengths of gaps between 1's of an
independent Bernoulli sequence arises. Pittel shows that the
sequence of lengths of macrosteps obtained when applying Wilson's
algorithm to the complete graph with $N$ vertices has the same
distribution as the sequence of spacings between successes of
independent Bernoulli($j/N$) variables, $2 \leq j \leq N$.
For each $N=1,2,\ldots$, let $(Y_{N,j}, j=1,2,\ldots)$
be such an independent sequence of Bernoulli random variables, where
for $1 \leq j \leq n$, $Y_{N,j}$ has parameter $j/N$, and for $j > N$,
$Y_{N,j}=1$. Let $Z_{N,1}$ be the smallest index $j$ such that $Y_{N,j} = 1$,
let $Z_{N,2}$ be the second-smallest index $j$ such that $Y_{N,j} = 1$,
and define $Z_{N,j}$ similarly. Then for any positive integer $m$, as
$N \to \infty$, $$\frac{1}{\sqrt{N}}(Z_{N,1}, \ldots, Z_{N,m}) \convd
(P_1, \ldots, P_m),$$ where $P_1 < P_2 < \cdots$ are the successive
points of an inhomogeneous Poisson point process on $[0,\infty)$
with rate $t$ at time $t$; that is, at every continuity point
of the distribution function of $(P_1, \ldots, P_m)$, the distribution
function of $(Z_{N,1}, \ldots, Z_{N,m})$ converges to the distribution
function of $(P_1, \ldots, P_m)$. This follows from, e.g., results of \cite[\S 2.6]{Dur}.

A similar type of scaling limit comes up
when analyzing repeat values. Let $R_{N,1} := R_N$ be the index of
first repeat in an i.i.d. sequence of random variables uniform on $[N]$,
and let $R_{N,2}$ be the second-smallest index $i$ such that
$X_i \in \{X_0, \ldots, X_{i-1}\}$, let $R_{N,3}$ be the third-smallest
such index, and so on. Then (see \cite{CP} and work cited there) the sequence $(R_{N,i}, i = 1,2,\ldots)$
has the same finite-dimensional scaling limits as the sequence $(Y_{N,j}, j = 1,2,\ldots)$:
for any positive integer $m$, as $N \to \infty$, $$\frac{1}{\sqrt N}(R_{N,1},
\ldots, R_{N,m}) \convd (P_1, \ldots, P_m).$$ This latter scaling limit
result is an integral part of one of Aldous's constructions
of the Brownian continuum random tree \cite{Al1}.

But now, using the fact that $\lambda_{N,k} \eqd (R_n - k|R_n > k)$
from Theorem \ref{seg}, we obtain the following:
\begin{Corollary} \label{lnkscale}
For a fixed $\mu > 0$, as $N \to \infty$, $(\lambda_{N,\lf\mu \sqrt
N\rf})/\sqrt N$ converges in distribution to $(P_1 - \mu|P_1 >
\mu)$, where $P_1$ is the first point of an inhomogeneous Poisson
point process on $[0,\infty)$ with rate $t$ at time $t$.
\end{Corollary}

This result can also be used to provide an estimate for $\nu_{N,k}$, the
number of colors in the loop-erased random walk from $N$ to $[k]$. It
is clear that conditional on $\lambda_{N,k} = j$, the expected
value of $\nu_{N,k}$ is given by $$E(\nu_{N,k}|\lambda_{N,k} = j) = 1 +
\sum_{i=1}^{j-1} \frac{1}{k+i+1}.$$ It thus follows that
\begin{eqnarray*}
E(\nu_{N,\lf\mu\sqrt N\rf}) &=& \sum_{j=1}^{N - \lf\mu\sqrt N\rf}
\left(1 + \frac{1}{\lf\mu \sqrt N\rf + 2} + \cdots + \frac{1}{\lf\mu
\sqrt N\rf + j}\right) \P(\lambda_{N,\lf\mu\sqrt N\rf} = j)
\end{eqnarray*}
But as $N \to \infty$, the term $1 + \frac{1}{\lf\mu \sqrt N\rf + 2}
+ \cdots + \frac{1}{\lf\mu \sqrt N\rf + j}$ approaches $1 +
\log\left(1 + \frac{j}{\mu \sqrt N}\right)$, and by Corollary
\ref{lnkscale}, the term $\P(\lambda_{N,\lf\mu\sqrt N\rf} = j)$
approaches $\frac{1}{\sqrt N} \left(\frac{j}{\sqrt N} + \mu\right)
\exp\left(-j^2 / 2N - \mu j/\sqrt N\right).$ Therefore, as $N \to
\infty$, \begin{eqnarray} \label{nulimit} E(\nu_{N,\lf\mu\sqrt
N\rf}) &\to& \int_0^\infty \left(1 + \log\left(
1+\frac{x}{\mu}\right) \right) (\mu + x)\exp\left(-\frac{x^2}{2} -
\mu x\right)
\\&=& 1 - \log \mu + \exp\left(\frac{\mu^2}{2}\right) \int_\mu^\infty t\log t \exp\left(
-\frac{t^2}{2}\right)dt
\end{eqnarray}
after some simplification and the substitution $t = x + \mu$.

An alternative derivation of this scaling limit comes from the fact
that the length of the loop-erased random walk increases at unit
speed until a length $j$ when a loop occurs, after which its new
length is uniformly distributed among $\{0, 1, \ldots, j\}$. As
such, it is closely related to the {\it standard Rayleigh process}
$(R_t, t \geq 0)$ \cite{EPW}, defined as follows: let $R_0 = 0$, and
for $P_{i-1} < t < P_i$ (with the convention that $P_0 = 0$), let
$R_t$ grow at unit speed; at each time $P_i$, let $R_{P_i}$ be
selected uniformly within the interval $(0, R_{P_i-})$. If we also
make note of the basic fact that $\z_{N,k}$, the geometric time at
which the walk from $N$ hits $[k]$, satisfies the scaling limit
$\z_{N,\lf\mu\sqrt N \rf} / \sqrt N \convd X_\mu$, where $X_\mu$ has
an exponential distribution with parameter $\mu$, then the following
is clear:

\begin{Corollary} \label{lnkscale2}
As $N \to \infty$, $\lambda_{N,\lf\mu\sqrt N\rf} / \sqrt N$
converges in distribution to $R_{X_\mu}$, where $(R_t, t \geq 0)$ is
the standard Rayleigh process, and $X_\mu$ is independent of $(R_t,
t \geq 0)$ and has an exponential distribution with parameter $\mu$.
\end{Corollary}

Using the fact from \cite{EPW} that $R_t \eqd P_1 \wedge t$, where
$P_1$ has the standard Rayleigh distribution, shows that the scaling
limits in Corollaries \ref{lnkscale} and \ref{lnkscale2} have the
same distribution.

Consider a finite $k$, condition on the length of the loop-erased
random walk from $N$ to $[k]$ equalling $j$, and let $C_1 < \cdots <
C_{\nu_{N,k}-1}$ be the color-changing indices, i.e., those indices
$i$ in $[j-1]$ such that $B_i = 1$. By Theorem \ref{seg}, it follows
that
\begin{eqnarray*} \P(C_1 > i) &=& \left(1 - \frac{1}{k+2}\right)
\left(1 - \frac{1}{k+3}\right) \cdots \left(1 - \frac{1}{k+i+1}\right) = \frac{k+1}{k+i+1}\\
\P(C_m > i+j|C_{m-1} = i) &=& \frac{k+i+1}{k+i+j+1}\end{eqnarray*}
If we now let $k = \lf\mu\sqrt N\rf$, and let $N \to \infty$, then
it follows that $C_1/\sqrt N \convd D_1$, where $P(D_1 > \lambda) =
\frac{\mu}{\mu + \lambda}$. More generally, we see that:

\begin{Proposition}
As $N \to \infty$, $$\frac{1}{\sqrt N}(C_1, \ldots,
C_{\nu_{N,\lf\mu\sqrt N\rf} - 1}) \convd (D_1, \ldots,
D_{\nu(\mu)}),$$ where $D_1 < D_2 < \cdots$ are the points of an
inhomogeneous Poisson point process with rate $\frac{1}{\mu + t}$ at
time $t$, and $\nu(\mu):= \sup\{i: D_i < R_{X_\mu}\}$.
\end{Proposition}

The sequence $(D_1, \ldots, D_{\nu(\mu)})$ is called the sequence of
{\it ladder indices} of the standard Rayleigh process up to time
$X_\mu$: if, at a jump time $P_i$, the Rayleigh process jumps down
to some $Q_i$ uniformly chosen in $(0, R_{P_i-})$, then the sequence
$(D_1, \ldots, D_{\nu(\mu)})$ is the subsequence $(Q_{i_1}, \ldots,
Q_{i_{\nu(\mu)}})$ of $(Q_1, Q_2, \ldots)$ up to time $X_\mu$ such
for each index $i_c$, $Q_{i_c} < Q_{i_c +m}$ for all $m =
1,2,\ldots$.

As a check, note that conditional on $R_{X_\mu} = x$, the number of
colors in the scaled walk is one plus a Poisson random variable with
parameter $\log \left(1 + x/\mu\right)$, and thus the expected value
of this number of colors is
$$1 + \int_0^\infty \log \left(1 + \frac{x}{\mu}\right) (x+\mu) \exp\left(-\frac{x^2}{2}
- \mu x\right) dx,$$ which is the same as derived in equation
(\ref{nulimit}).

\subsection{Stick-breaking}
In the previous subsection, letting $k = \lf \mu \sqrt N \rf$ led to
scaling limit results which were closely related to particular
Poisson and Rayleigh processes on $[0,\infty)$. Fixing $k$ while
letting $N \to \infty$ leads to a very different result. To explore
this, first note that Theorem \ref{seg} also relates to some other
basic results on random compositions. See \cite[\S I.5]{Fel} and
\cite[pp. 52-53]{Dur} for discussions of how the record times of a
sequence of i.i.d. random variables $(W_i, i= 1, 2, \ldots)$---i.e.,
the times $j$ such that $W_j > W_i$ for all $1 \leq i \leq j-1$---is
distributed like the occurrences of 1's in a sequence of independent
Bernoulli$(\frac{1}{i+1})$ random variables.

In the setting of Theorem \ref{seg}, let $(\s_{N,k,1}, \ldots,
\s_{N,k,\nu_{N,k}})$ denote the random composition of
$\lambda_{N,k}$ representing the numbers of edges of each color,
working along the colored loop-erased path from $N$ to $[k]$, and
consider the reversed sequence of segment lengths \be (\s'_{N,k,1},
\ldots, \s'_{N,k,\nu_{N,k}}) := (\s_{N,k,\nu_{N,k}}, \ldots,
\s_{N,k,1}). \ee Also, suppose that we condition on $\lambda_{N,k} =
j$ for some $1 \leq j \leq N-k$, so that if $C_i$ for $1 \leq i \leq
j-1$ is the indicator of whether $i$ is a partial sum of the
sequence $(\s'_{N,k,m}, 1 \leq m \leq \nu_{N,k})$, then the
variables $(C_i, 1 \leq i \leq j-1)$, conditioned on $\lambda_{N,k}
= j$, are independent Bernoulli$(\frac{1}{k+j+1-i})$ random
variables. We know from the proof of Theorem \ref{seg} that
$\P(\s'_{N,k,1} = j) = \frac{k+1}{k+j}$ (since this corresponds to
there being only one color in the loop-erased walk); however,
because of Theorem \ref{seg}, we now also see that for $1 \leq n
\leq j-1$, $$\P(\s'_{N,k,1} = n) = \left(1 -
\frac{1}{k+j}\right)\left(1 - \frac{1}{k+j-1}\right) \cdots \left(1
- \frac{1}{k+j+2-n}\right) \frac{1}{k+j+1-n} =\frac{1}{k+j}$$ so
that $\s'_{N,k,1}$ is equally likely to be any of $1,2,\ldots, j-1$.

It is easy to generalize this, to see that conditional on
$(\s'_{N,k,1}, \ldots, \s'_{N,k,m}) = (k_1,\ldots, k_m)$ with $c:= k_1
+ \cdots +k_m < j$, \ba \P(\s'_{N,k,m+1} = j-c) &=&\frac{k+1}{k+j-c}\\
\P(\s'_{N,k,m+1} = k_{m+1})&=& \frac{1}{k+j-c}, ~~1 \leq k_{m+1}
\leq j-c-1\ea Now since $\lambda_{N,k} \to \infty$ in probability as $N
\to \infty$, we see that as $N \to \infty$, \be \label{stick}
\lambda_{N,k}^{-1}(\s'_{N,k,1}, \ldots, \s'_{N,k,\nu_{N,k}}) \convd
(U_1,(1-U_1)U_2, (1-U_1)(1-U_2)U_3, \ldots) \ee where the $(U_i, i
= 1,2,\ldots)$ are independent uniform$(0,1)$ random variables.
The right hand side of (\ref{stick}) is known as the {\it
continuous uniform stickbreaking process} defined by the
$(U_i,i=1,2,\ldots)$.

\subsection{Markovian and non-Markovian properties of the colored walk}
Although the sequence of lengths $L_n$ of the loop-erased walk is a Markov chain,
we will now argue that the sequence of compositions of $L_n$ defined by the lengths of colored segments of the loop-erased
walk is not a Markov chain.

Recall that a {\em composition} of a positive integer $\ell$ is a
sequence of positive integers with sum $\ell$. Here the terms of the
sequence represent the lengths of stretches of edges of the same
color in a loop-erased walk with total length $\ell$. Such a
composition of $\ell$ is conveniently encoded by the string of
$\ell$ binary {\em bits} $(y_1, \ldots, y_\ell)$ where $y_1 = 1$ and
$y_i$ is the indicator of a color change between the $(i-1)^{th}$
and $i^{th}$ edges. So the number of compositions of $\ell$ is
$2^{\ell-1}$. If $y = (y_1, \ldots, y_\ell)$ and $z = (z_1, \ldots,
z_m)$ are compositions of $\ell$ and $m$ respectively, call $y$ a
{\em truncation} of $z$ if $\ell \le m$, $y_i = z_i $ for $1 \le i
\le \ell-1$, and $y_\ell \leq z_\ell$. We introduce also a {\em
trivial composition}, corresponding to a sequence with no terms,
which is regarded as a truncation of every composition of a positive
integer.

For $n = 0,1,\ldots$, let $C_n$ denote
the composition induced by coloring the segments of the loop-erasure
of $(X_0, \ldots, X_n)$, where it is understood that if the loop-erasure
has one vertex and no edges, then $C_n$ is the trivial composition.
Observe that $X_{n+1}$ belongs to the set of values in the loop-erasure of $(X_0, \ldots, X_n)$ if and only if
$C_{n+1}$ is a truncation of $C_n$. The sequence $(C_n, n = 0,1,\ldots)$ has the following dynamics:
\begin{itemize}
\item If $C_n$ is the trivial composition, then $C_{n+1}$ is
the trivial composition with probability $1/N$, and $C_{n+1} = (1)$
with probability $1 - 1/N$.
\item If $C_n$ is some non-trivial composition $c$ of $\ell$,
then $C_n$ is either some truncation of $C_{n-1}$ or an extension of $C_{n-1}$ by one term,
according to whether or not $X_n$ creates a loop:
\begin{itemize}
\item if $X_n$ creates a loop, then $C_{n+1}$ extends $C_n$ by adding the bit $1$ with probability $(N-\ell-1)/{N}$,
while $C_{n+1}$ is equally likely to be each of the $\ell +1$ possible truncations of $C_n$;
\item if $X_n$ does not create a loop, then
$C_{n+1}$ extends $C_n$ by adding the bit $0$ with probability $(N-\ell-1)/{N}$,
while $C_{n+1}$ is equally likely to be each of the $\ell +1$ possible truncations of $C_n$;
\end{itemize}
\end{itemize}
The sequence $(C_n, n = 0,1,\ldots)$ is not
a Markov chain because
the transition dynamics from $C_n$ to $C_{n+1}$ depend on whether $X_n$ created a loop, which is determined
by the relationship between $C_{n-1}$ and $C_n$. However,
this analysis does show that the sequence of pairs of compositions,
$((C_n,C_{n+1}), n = 0,1,\ldots)$ is a Markov chain.

\section{Open questions} \label{open}
This article has examined loop-erased random walk which stops when
we reach some marked subset of $[N]$ of size $k$, at which point the walk has a
random length $\lambda_{N,k}$; in this case, we showed that
conditional on the value of $\lambda_{N,k}$, the composition of
color segments was distributed as the length of spacings between 1's
of a sequence of Bernoulli random variables. What if we stop the
walk when the length of the loop-erasure reaches $m$, for some $1
\leq m \leq N-k$? Will the resulting composition of colors have a
similar distribution, as the lengths of spacings? How will the
composition be distributed if we just stop
at some fixed finite time $m$, and look at the colored loop-erasure
of $(X_0,\ldots, X_m)$? The situation appears to be like the Ray-Knight description
of the distribution of Brownian local times, where results are much simpler
for suitable stopping times than for fixed times.

In Section \ref{scalims}, we derived a number of results which hold
at the end of the stopped walk. But in this situation too, what
happens at an intermediate stage? In the case of $k = \lf\mu\sqrt
N\rf$, it seems that the length of the loop-erased random walk,
considered as a stochastic process, should converge to the standard
Rayleigh process.

Finally, it seems that there is some sort of ``critical'' behavior
when $k$ is of order of magnitude $\sqrt N$, as detailed in Section
\ref{scalims}. What happens to the colored walk when $k = o(\sqrt
N)$ or $k = \omega(\sqrt N)$? In the former case, the walk does not
get stopped by an exponential time, and it appears that the number
of colors of the stopped walk should increase as $\log N$. In the
latter case, in some sense the walk gets stopped before it can make
loops, and it appears that the number of colors should converge to 1
almost surely.

\end{document}